\documentstyle[12pt]{article} 
\newfont{\twlvmsb}{msbm10 scaled\magstep1}
\newfont{\ninemsb}{msbm9}
\newfont{\sixmsb}{msbm6}
\newfam\msbfam
\textfont\msbfam=\twlvmsb
\scriptfont\msbfam=\ninemsb
\scriptscriptfont\msbfam=\sixmsb
\catcode`\@=11
\def\Bbb{\ifmmode\let\next\Bbb@\else
  \def\next{\errmessage{Use \string\Bbb\space only in math mode}}\fi\next}
\def\Bbb@#1{{\Bbb@@{#1}}}
\def\Bbb@@#1{\fam\msbfam#1}
\newfont{\largeeufm}{eufm10 scaled\magstep4}
\newfont{\twlveufm}{eufm10 scaled\magstep1}
\newfont{\elveufm}{eufm10 at 11pt}
\newfont{\teneufm}{eufm10}
\newfont{\nineeufm}{eufm9}
\newfam\eufam
\textfont\eufam=\twlveufm
\scriptfont\eufam=\nineeufm
\def\frak{\ifmmode\let\next\frak@\else
\def\next{\errmessage{Use \string\frak\space only in math mode}}\fi\next}
\def\frak@#1{{\fam\eufam{{#1}}}}
\catcode`@=12
\newcommand{\Z}{{\Bbb Z}} 
\newcommand{\Zz}{{\Bbb Z}_2} 
\newcommand{\C}{{\Bbb C}} 
\newcommand{\gl}{{\frak g}} 
\newcommand{\Uq}{U_q({\frak g})} 
\newcommand{\U}{U({\frak g})} 
\newcommand{\Uqp}{U_q({\frak p})} 
\newcommand{\Upm}{U_q({\frak p}_\pm)} 
\newcommand{\Up}{U_q({\frak p}_+)} 
\newcommand{\Um}{U_q({\frak p}_-)} 
\newcommand{\Ul}{U_q({\frak l})} 
\newcommand{\Sq}{{\Bbb S}_q^{m|n-1}} 
\newcommand{\CP}{{\Bbb C}P_q^{m|n-1}} 
\newcommand{\ba}{\begin{eqnarray}}
\newcommand{\na}{\end{eqnarray}}
\newcommand{\ban}{\begin{eqnarray*}}
\newcommand{\nan}{\end{eqnarray*}}
\newtheorem{lemma}{Lemma}
\newtheorem{proposition}{Proposition}

\begin{document} 
\title{\normalsize{\bf STRUCTURE  AND REPRESENTATIONS OF\\
 THE QUANTUM   GENERAL   LINEAR SUPERGROUP} } 
\author{ \small R. B. ZHANG\\ 
\small Department of Pure Mathematics\\
\small  University of Adelaide\\
\small  Adelaide, Australia}  
\date{} 
\maketitle
\begin{abstract} 
The structure and representations of the quantum general linear 
supergroup $GL_q(m|n)$ are studied systematically 
by investigating the Hopf superalgebra $G_q$ of its representative 
functions. $G_q$ is factorized into $G_q^{\pi}\ G_q^{\bar\pi}$, 
and a Peter - Weyl basis is constructed for each factor. 
Parabolic induction for the quantum supergroup is developed.  
The underlying geometry of induced representations    
is discussed, and an analog of Frobenius reciprocity is obtained. 
A quantum Borel - Weil theorem is proven
for the covariant and contravariant tensorial irreps, and     
explicit realizations are given  for classes of
tensorial irreps in terms of sections of quantum
super vector bundles over quantum projective superspaces.
\end{abstract} 

\vspace{1cm}
\section{\normalsize INTRODUCTION} 
Quantized universal enveloping superalgebras \cite{Bracken}
\cite{Chaichian} ( which will be called quantum superalgebras for simplicity) 
represent the most important generalizations of the Drinfeld - Jimbo
\cite{Drinfeld}   
quantized universal enveloping algebras.
Their origin can be traced back to the Perk - Schultz solution 
of the Yang - Baxter equation and also the work of 
Bazhanov and  Shadrikov\cite{Perk}.   
However, systematical investigations of such algebraic structures 
only started about six years ago, but in an intensive manner.  
By now the subject has been developed quite extensively: 
the quasi - triangular Hopf superalgebraic 
structure of the quantum superalgebras was investigated \cite{Double}; 
the representation theory of large classes of quantum 
(affine)superalgebras and super Yangians was developed 
\cite{I}\cite{II};  
applications of quantum superalgebras to integrable  
two dimensional models in statistical mechanics and quantum 
field theory were extensively explored \cite{Bracken}\cite{Delius}. 
Quantum superalgebras have also been applied to the study of 
knot theory and $3$ - manifolds \cite{Links}\cite{Manifolds}, 
yielding many new topological invariants, notably, the 
multi - parameter generalizations of Alexander - Conway polynomials.

Closely related to the Drinfeld - Jimbo algebras are the quantum 
groups introduced by Woronowicz 
and Faddeev - Reshetikhin - Takhatajan \cite{Faddeev}, 
which are, in the spirit 
of Tannaka - Krein duality theory, the `groups' associated 
with the quantized universal enveloping algebras. 
One very important aspect of quantum groups is their 
geometrical significance: they provide a concrete framework 
for developing noncommutative geometry \cite{Manin},  
in particular, for investigating notions such as quantum 
flag varieties \cite{Lakshmibai} and quantum fibre bundles.

Our aim  here is to study the structure and representations 
of the  quantum  general  linear supergroup $GL_q(m|n)$ in a   
systematical fashion by investigating the algebra of its 
representative functions.  We start in section 2 with a concise 
treatment of finite dimensional unitary representations of 
$U_q(gl(m|n))$. Results will be repeatedly used in the remainder of the 
paper.  In section 3 we define the the  quantum  general  linear
supergroup $GL_q(m|n)$, or more exactly, the superalgebra $G_q$ 
of functions on it. This is done by first defining the bi - 
superalgebras $G_q^\pi$ and $G_q^{\bar\pi}$, which are respectively 
generated by the matrix elements of the vector irrep and 
its dual irrep.   Peter - Weyl type of bases for these 
bi - superalgebras are constructed. The $G_q$ is defined to 
be generated by $G_q^\pi$ and $G_q^{\bar\pi}$ with some 
extra  relations. It has the structures of a $\ast$ - Hopf superalgebra, 
which separates points of $U_q(gl(m|n))$,  and factorizes 
into $G_q^\pi G_q^{\bar\pi}$.  Section 4 treats the representation 
theory of the quantum supergroup, and in particular, parabolic 
induction. The geometrical interpretation of induced representations 
is discussed,  leading  naturally to the concepts of quantum 
homogeneous spaces and quantum super vector bundles.  
A quantum analog of Frobenius reciprocity is obtained; 
and a quantum version of the Borel - Weil theorem is proven 
for the covariant and contravariant tensorial irreps. Section 5 
gives the explicit realizations of two infinite classes of 
tensorial irreps in terms of sections of quantum 
super vector bundles over the quantum projective superspace. 
In doing this, we also treat the quantum projective superspace 
in some detail.

\section{\normalsize UNITARITY REPRESENTATIONS OF $U_q(gl(m|n))$}  
The finite dimensional unitary representations of $U_q(gl(m|n))$ 
were classified in \cite{Scheunert}.  
Here we will reformulate the results on the covariant and 
contravariant tensor irreps  so that they can be 
readily used in the remainder of the paper. 
The material presented here also heavily relies on 
references \cite{I} and \cite{Gould}.  

\subsection{\normalsize   Hopf $\ast$ - superalgebras 
and unitary representations}
Let $A$ be a $\Zz$ - graded associative algebra over the complex field 
$\C$.   Its underlying $\Zz$ - graded vector space is the 
direct sum $A=A_0\oplus A_1$ of the even subspace 
$A_0$ and the odd subspace $A_1$. 
We introduce the grading index $[\ ]: A_0\cup A_1\rightarrow \Zz$  
such that $[a]=\theta$ if $a\in A_\theta$.  We will call $A$ a 
$\Zz$ - graded $\ast$ - algebra, or $\ast$ - superalgebra,  
if there exists 
an even anti - linear anti - automorphism $\ast: A\rightarrow A$ 
such that $\ast\circ \ast = id_A$.  
We will denote $\ast(a)$ by $a^\ast$. 
Needless to say, $\ast(a b)=b^\ast a^\ast$, $a, b\in A$. 

An important new feature of the $\Zz$ - graded case is that 
for a given $*$ - operation of $A$, there exists an associated 
$*'$ such that 
\ba 
*'(a)&=(-1)^{[a]} a^*, \label{II} 
\na     
for $a$  being homogeneous,  
and extends to the whole of $A$ anti - linearly.   
There also exist the so called graded $*$ - operations, which, however, 
are not useful for this paper, thus will not be discussed any further.

Let $A$ and $B$ be two $\Zz$ - graded $\ast$ - algebras. 
Then $A\otimes_\C B$ has a natural $\Zz$ - graded $\ast$ - algebra 
structure, with the $\ast$ - operation defined for homogeneous 
elements by  
\ban 
\ast(a\otimes b) &=& (-1)^{[a][b]} a^\ast\otimes b^\ast, 
\nan 
and for all the elements by extending this anti - linearly. 

Consider a $\Zz$ - graded Hopf algebra ( also called Hopf superalgebra ) 
H, with multiplication $m$, unit $1_H$, co - multiplication 
$\Delta$, co - unit $\epsilon$ and antipode $S$. We emphasize that 
the antipode is a {\em linear} anti- automorphism of the underlying 
algebra of $H$. In particular, for homogeneous $a, b\in A$, 
we have $S(a b) = (-1)^{[a][b]}S(b) S(a)$.  $H$ will be called a 
$\Zz$ - graded Hopf $\ast$ - algebra, or Hopf $\ast$ - superalgebra, 
if the underlying algebra of $H$ is a $\ast$ - superalgebra such that 
$\Delta$ and $\epsilon$ are $\ast$ - homomorphisms, i.e., 
\ban 
\ast\circ \Delta = \Delta\circ \ast, & 
\ast\circ\epsilon=\epsilon \circ \ast. 
\nan   
These properties together with the defining relations of the antipode  
\ban 
m\circ ( S\otimes id ) \Delta &  
= m\circ (id\otimes S) \Delta  &= 1_H\epsilon 
\nan 
imply that  
\ban 
S\circ\ast\circ S\circ\ast &=& id_H. 
\nan 

Let $V$ be a left $H$ - module. If there exists a non - degenerate 
sequilinear form $(\ , \ ):$  $V\otimes V\rightarrow \C$, such that 
\ban 
(i).& (a v, \ u)=(v,\ a^* u), &\forall u, v\in V, \ a\in H,\\  
(ii).&  (v,\ v)\ge 0 , &(v,\ v)=0 \ \  \mbox{iff}\  \ v=0,  
\nan  
we call $V$ and the associated representation of $H$ unitary.

Unitary representations have the following important properties\\ 
{\em i).   A unitary representation is completely reducible;} \\  
{\em ii).  The tensor product of two unitary 
       ( with respect to the same $*$ - operation ) 
        representations  is again unitary;}\\   
{\em iii). If a representation is unitary with respect to $*$, then its dual
       is unitary with respect to $*'$.}

All the three assertions are well known, but there are some related 
matters worth discussing.  One is concerned with  
the requirement that two representations must be unitary with 
respect to the same $*$ - operation in order for their tensor 
product to be unitary as well.  The tensor product $V\otimes_\C W$ 
of two $H$ - modules has a natural $H$ module structure 
\ban 
a \{v\otimes w\} &=& \sum_{(a)} (-1)^{[a_{(2)}] [v]} 
                     a_{(1)} v \otimes a_{(2)} w.   
\nan 
If both $V$ and $W$ are equipped with sequilinear forms 
$( \ , \ ): V\otimes_\C V \rightarrow \C$,  
and $( \ , \ ): W\otimes_\C W \rightarrow \C$,  
we can define a sequilinear form $(( \ , \ )): (V\otimes_\C W)^{\otimes 2}
               \rightarrow \C$ by  
\ban 
((v_1\otimes w_1,\ v_2\otimes w_2)) &=& (v_1,\ v_2) (w_1,\ w_2).  
\nan 
Now if both $V$ and $W$ are unitary with respect to the same $*$ - 
operation, then $(( \ , \ ))$ is clearly positive definite and 
nondegenerate. Furthermore,    
\ban ((v_1\otimes w_1,\ a\{v_2\otimes w_2\})) 
&=& \sum_{(a)} (-1)^{[a_{(2)}][v_2]} 
(( a_{(1)}^* v_1\otimes a_{(2)}^* w_1, \ v_2\otimes w_2)) \\ 
&=& (( a^*\{v_1\otimes w_1\}, \ v_2\otimes w_2)). 
\nan   
Therefore, $V\otimes_\C W$ indeed furnishes a unitary $H$ - module. 
On the other hand, if, say, $V$ is $*$ - unitary, while $W$ is 
$*'$ - unitary, then one can easily see that the above calculations 
will fail to go through.   

The other concerns the third assertion, the validity of which 
actually requires some qualification, namely, 
the Hopf * - superalgebra $H$ in question must admit an even group like 
element $K_{2\rho}$ satisfying  
\ba 
K_{2\rho}^* = K_{2\rho}, &  
S^2(a) = K_{2\rho} a K_{2\rho}^{-1}, \ \ \ \forall a\in H. 
\label{antipode}  
\na   
Let $V$ be a locally finite module over $H$, which is unitary 
with respect to the sequilinear form $(\ , \ ): V\otimes V 
\rightarrow \C$.  For every $v\in V$, we define $v^\dagger$ by 
$v^\dagger ( w ) = (v, w)$, $\forall w\in V$, and denote the linear 
span of all such $v^\dagger$ by  $V^\dagger$, which is a subspace 
of the dual vector space of $V$.  The $V^\dagger$ has a natural 
$H$ module structure, with the action of $H$ given by 
\ban 
(a v^\dagger)(w) &=& (-1)^{[a][v^\dagger]} v^\dagger( S(a) w),   
\ \ \ w\in V. 
\nan 
Unitarity of $V$ leads to 
\ban 
a v^\dagger &=& (-1)^{[a][v]}\left( * S (a) v \right)^\dagger. 
\nan 
We define a sequilinear form
$( \ , \ )^\prime: V^\dagger\otimes V^\dagger\rightarrow \C$ by 
\ban 
(v^\dagger, \ w^\dagger)^\prime &=& (K_{2\rho} w, \  v). 
\nan    
It follows from the properties of the original form on $V$ that 
$(\ , \ )^\prime$ is positive definite and nondegenerate.  
A straightforward calculation shows that  
\ban 
(a v^\dagger, \ w^\dagger)^\prime&=& ( v^\dagger, *'(a) w^\dagger )^\prime, 
\nan 
where $*'$ is defined by (\ref{II}).

\subsection{\normalsize  $U_q(gl(m|n))$} 
Throughout the paper, we will denote by $\frak g$ the complex 
Lie superalgebra $gl(m|n)$, and by $U({\frak g})$ 
its universal enveloping algebra. 
As is well known, there are the Drinfeld and Jimbo two versions 
of the  quantized universal enveloping algebra $\Uq$ of 
$\frak g$, which, though, have very similar properties at generic $q$.  

It is the Jimbo version of $\Uq$ that will be used in this paper. 
Now $\Uq$ is a $\Zz$ - graded unital associative algebra 
over $\C(q, q^{-1})$, $q$ being an indeterminate, 
generated by $\{K_a,  \ K_a^{-1}, \ a\in {\bf I};   \ 
E_{b\   {b+1}},$ $ \ E_{b+1,   b}, \ b\in {\bf I}'\}$, 
${\bf I}=\{1, 2, ..., m+n\}$, ${\bf I}'=\{1, 2, ..., m+n-1\}$,   
subject to the following relations  
\ba 
K_a K_a^{-1}=1, 
& & K_a^{\pm  1} K_b^{\pm 1} = K_b^{\pm 1}  K_a^{\pm 1}, \nonumber \\ 
K_a E_{b\ b\pm 1} K_a^{-1} &=& 
q_a^{\delta_{a b} -\delta_{a\  b\pm 1}} E_{b\  b\pm 1}, \nonumber \\ 
{}[E_{a\,  a+1},\,E_{b+1\,  b}\}& =& \delta_{a b}
(K_a K_{a+1}^{-1} - K_a^{-1} K_{a+1})/(q_a - q_a^{-1}),\nonumber \\ 
(E_{m\, m+1})^2 &=& (E_{m+1\, m})^2 = 0, \nonumber \\  
E_{a\,  a+1} E_{b\,  b+1} &=& E_{b\,  b+1} E_{a\,  a+1},\nonumber \\     
E_{a+1\, a} E_{b+1\, b} &=&E_{b+1\, b} E_{a+1\, a}, \ \ \  
\vert a - b\vert \ge 2, \nonumber \\ 
{\cal S}^{(+)}_{a \ a\pm 1}&=&{\cal S}^{(-)}_{a \ a\pm 1}=0,  
\ \ \ a\ne m,\nonumber \\ 
\{ E_{m-1\, m+2},\ E_{m\, m+1}\} &=& 
\{ E_{m+2\, m-1},\ E_{m+1\, m}\} = 0,  \label{quantum}  
\na 
where $q_a=q^{ (-1)^{[a]} }$,  
\ban 
{\cal S}^{(+)}_{a \ a\pm 1}&=& 
(E_{a\, a+1})^2  E_{a\pm 1\, a+1\pm 1} - (q +
q^{-1}) E_{a\, a+1} \ E_{a\pm 1\, a+1\pm 1} \ E_{a\, a+1}\\ 
& +& E_{a\pm 1\, a+1\pm1 }\ (E_{a\, a+1})^2,    \\ 
{\cal S}^{(-)}_{a \ a\pm 1}&=& 
(E_{a+1\, a})^2\,E_{a+1\pm 1\, a\pm 1} - (q +
q^{-1}) E_{a+1\, a}\ E_{a+1\pm 1\, a\pm 1} \ E_{a+1\, a}\\ 
&+& E_{a+1\pm 1\, a\pm 1}\ (E_{a+1\, a})^2, 
\nan 
and $E_{m-1\, m+2}$ and $E_{m+2\, m-1}$ are the $a=m-1$, $b=m+1$, 
cases of the following elements 
\ban 
E_{a\, b} &=& E_{a\, c} E_{c\, b} - q_c^{-1} E_{c\, b} E_{a\, c}, \\  
E_{b\, a} &=& E_{b\, c} E_{c\, a} -    q_c   E_{c\, a} E_{b\, c}, 
\ \ \ a<c<b. 
\nan 
The $\Zz$ grading of the algebra is specified such that 
the elements $K_a^{\pm 1}$, $\forall a\in {\bf I}$, 
and $E_{b\,  b+1}$, $E_{b+1\,  b}$,  $b\ne m$, are even, 
while $E_{m\, m+1}$ and $E_{m+1\, m}$ are odd. 
Above, we have also used the 
notation $[a]=\left\{\begin{array}{l l} 
               0,  & \mbox{if} \ a\le m, \\ 
               1,  & \mbox{if} \ a>m. 
              \end{array}\right. $ 

On the other hand, the Drinfeld version of $\Uq$ is defined over 
$\C[[\hbar]]$, $q=\exp(\hbar)$, and is completed with respect 
to the $\hbar$ - adic topology of $\C[[\hbar]]$. It is generated 
$\{E_{a\, a}, \ a\in {\bf I};   \
E_{b\   {b+1}},$ $ \ E_{b+1,   b}, \ b\in {\bf I}'\}$,
subject to the same relations (\ref{quantum}) with 
$$K_a = q_a^{E_{a\, a}}.$$

It is well known that $\Uq$ has the structure of a $\Zz$ graded 
Hopf algebra, with a co - multiplication 
\ban 
\Delta(E_{a\, a+1}) &=& E_{a\,  a+1} \otimes
K_a K_{a+1}^{-1} + 1 \otimes E_{a\, a+1}, \\ 
\Delta(E_{a+1\, a}) &=& E_{a+1\, a }\otimes 1 + K_a^{-1} K_{a+1}  
\otimes E_{a+1\, a}, \\
\Delta(K_a^{\pm 1}) &=&K_a^{\pm 1}\otimes K_a^{\pm 1},
\nan 
co - unit
\ban 
\epsilon(E_{a\, a+1})&=&E_{a+1\, a}=0, \ \ \forall a\in{\bf I}', \\ 
\epsilon(K_b^{\pm 1})&=&1,  \ \ \ \forall b\in{\bf I}, 
\nan  
and antipode
\ban 
S(E_{a\, a+1}) &=& - E_{a\, a+1} K_a^{-1} K_{a+1}, \\ 
S(E_{a+1\, a}) &=& - K_a K_{a+1}^{-1}E_{a+1\, a}, \\  
S(K_a^{\pm 1}) &=&K_a^{\mp 1}\otimes K_a^{\mp 1}. 
\nan

At generic $q$, the Jimbo version of $\Uq$ has more or less the 
same representation theory as that of the Drinfeld version \cite{I}.  
Let $\{\epsilon_a | a\in{\bf I}\}$ be the basis of a vector space
with a bilinear for
$  (\epsilon_a,\ \epsilon_b )=(-1)^{[a]}\delta_{a b}$.
The roots of the classical Lie superalgebra $gl(m|n)$ can be
expressed as
\ban
\epsilon_a - \epsilon_b, & a\ne b,& a, \, b\in{\bf I}.
\nan
For later use, we define
\ban
2\rho&=&\sum_{a\le b} (-1)^{[a]+[b]} ( \epsilon_a - \epsilon_b ).
\nan
 
From \cite{I} we know that
every finite dimensional irreducible $\Uq$ module is of highest weight
      type and is essentially uniquely characterized by a
      highest weight. 
Let $W(\lambda)$ be an irreducible $\Uq$ module
with highest weight 
$\lambda =\sum_{a}\lambda_a \epsilon_a$, $\lambda_a\in\C$.  
There exists a unique ( up to
scalar multiples) vector $v^\lambda_+\ne 0$ in $W(\lambda)$, called
the highest weight vector, such that
  \ban E_{a a+1} v^\lambda_+&=& 0, \ \ \ a\in{\bf I}',  \\
       K_b v^\lambda_+&=& q_b^{\lambda_{b}} v^\lambda_+,  \ \ \ b\in{\bf I}.
   \nan
$W(\lambda)$   is finite dimensional if and only if $\lambda$ satisfies
   $\lambda_a - \lambda_{a+1}\in\Z_+$,  $a\ne m$, and in that case,
it has the same weight space decomposition as that of the corresponding
irreducible $gl(m|n)$ module with the same highest weight.

\subsection{\normalsize Unitarity of covariant and contravariant 
           tensor irreps}
From this section on, we will assume that $\Uq$ is obtained 
from the the Jimbo algebra by specializing $q$ to a {\em real positive 
parameter} different from $1$.    
To construct a $*$ - operation for $\Uq$, we first consider the 
Hopf subalgebra generated by $e=E_{a \ a+1}$, $f=E_{a+1 \ a}$, 
and $k=K_a K_{a+1}^{-1}$, for a fixed $a\ne m$. It is not difficult to 
show that $*(e)=f k$, $*(f)=k^{-1} e$, $*(k^{\pm 1})=k^{\pm 1}$ 
defines a $*$ - operation for this $U_q(sl(2))$ subalgebra. 
Possible generalizations of this to $\Uq$ are 
\ba 
{*(E_{a\, a+1})}&=& (-1)^{(\theta+1)\delta_{m a}} 
 E_{a+1\, a} K_a K_{a+1}^{-1}, \nonumber\\ 
{*(E_{a+1\, a})}&=& (-1)^{(\theta+1)\delta_{m a}} 
         K_a^{-1} K_{a+1} E_{a\, a+1}, \nonumber\\ 
{*(K_a^{\pm 1})}&=& K_a^{\pm 1},   
\na 
where $\theta = 1 $ or $2$.  
It is quite obvious that the `quadratic' relations of 
(\ref{quantum}) are preserved by the $*$ - operations, and   
we have also explicitly checked that the `Serre relations' 
are  preserved as well.  We will call the $*$ - operations type
1 and type 2 respectively when  $\theta=1$ and $2$.   

It is also well known that  
\ban 
K_{2\rho} &=&\prod_{a<b} \left( K_a \, K_b^{-1}\right )^{(-1)^{[a]+[b]}} 
\nan 
satisfies equation (2).

Now we consider the covariant and contravariant tensor irreps of 
$\Uq$. 
The vector irrep $\pi$ of $\Uq$ is of highest weight $\epsilon_1$. 
The corresponding module $\Bbb E$  has the 
standard basis $\{v_a | a\in{\bf I}\}$, such that 
\ban 
K_a v_b &=& q_a^{\delta_{a b}} v_b, \\ 
E_{a \ a\pm 1} v_b &=& \delta_{b\ a\pm 1} v_a. 
\nan 
Define a sequilinear form on ${\Bbb E}\otimes {\Bbb E}$ by  
\ban 
(v_a, \ v_b )&=& \delta_{a b} \prod_{c=1}^{a-1} q_c^{-1}. 
\nan  
Then it is  straightforward to show that with respect to the type 1 
$*$ - operation, we have 
\ban 
(E_{a \ a\pm 1} v_b,\ v_c ) &=& ( v_b,\ E_{a \ a\pm 1}^* v_c ), \\
(K_a v_b,\ v_c ) &=& ( v_b,\ K_a v_c ). 
\nan 
Therefore, the vector irrep is unitary of type 1. 

The $\Uq$ modules  ${\Bbb E}^{\otimes k}$, $k\in\Z_+$ ( ${\Bbb E}^0= \C$ ),  
obtained by repeated tensor products of the vector module with 
itself can be decomposed into  direct sums of irreducible 
type 1 unitary modules, and we will call each direct summand  
an irreducible contravariant tensor module, 
and the corresponding irreducible representation a 
contravariant tensor irrep. 

The contravariant tensor irreps can be characterized in the following way. 
Let $\Z_+$ be the set of nonnegative integers. 
Define a subset $\cal P$ of ${\Z_+}^{\otimes(m+n)}$ by 
\ban 
{\cal P}&=&\{ p=(p_1, p_2, ..., p_{m+n})\in {\Z_+}^{\otimes(m+n)} 
\mid\ p_{m+1}\le n, \  p_a\ge p_{a+1},\   a\in{\bf I}'\}.   
\nan  
We associate with each $p\in{\cal P}$, 
a $\lambda^{(p)}=\sum_{a =1}^{m+n} \lambda_a \epsilon_a$ defined by 
\ban  
\lambda_a=p_a, & &a\le m, \nonumber\\
\sum_{\mu=1}^n \lambda_{m+\mu} \epsilon_{m+\mu} 
&=& \sum_{\nu=1}^n \sum_{\mu=1}^{p_{m+\nu}} \epsilon_{m+\mu}.   
\nan   
Introduce the set 
\ba 
\Lambda^{(1)}&=&\{ \lambda^{(p)} \mid  p\in {\cal P} \}. 
\na 
{}From results of\cite{Scheunert} \cite{Gould}  \cite{I} we know
that an irrep of $\Uq$ is a contravariant  tensor irrep if and only if 
its highest weight belongs to $\Lambda^{(1)}$. 
Needless to say, all such irreps are type 1 unitary.

Let  $W(\lambda)$  be  an irreducible contravariant tensor  
$\Uq$ module with highest weight $\lambda\in\Lambda^{(1)}$. 
We define   $\bar{\lambda}$ to be its lowest weight, and 
set $\lambda^\dagger=-\bar{\lambda}$. An explicit formula 
for $\lambda^\dagger$ was given in \cite{Gould}
( section III. B. ), where a more compact characterization
for the sets $\Lambda^{(1)}$ and $\Lambda^{(2)}$, which will be 
defined presently,  was also given.  
We refer to that paper for details.
Now  the dual module $W(\lambda)^\dagger$ of $W(\lambda)$, 
which will we call a covariant tensor 
module, has highest weight $\lambda^\dagger$. 

We introduce the set 
\ba 
\Lambda^{(2)}&=& \{\lambda^\dagger\ | \ \lambda\in\Lambda^{(1)}\}.    
\na 
An irrep is a covariant tensor irrep if and 
only if its highest weight is contained in $\Lambda^{(2)}$, and 
all covariant tensor irreps are unitary of type 2.  
The most important example of type 2 unitary modules is 
the covariant vector module ${\Bbb E}^\dagger$, which 
is the dual of the vector module ${\Bbb E}$. Its  
highest weight is given by $-\epsilon_{m+n}$.

We summarize our discussions into the following 
\begin{proposition}\label{tensor}    
\begin{enumerate}
\item Each $\Uq$ module ${\Bbb E}^{\otimes k}$ 
    ( resp. $({\Bbb E}^\dagger)^{\otimes k}$ ), 
   $k\in\Z_+$, can be decomposed into a direct sum of irreducible 
   modules with highest weights belonging to $\Lambda^{(1)}$ 
   ( resp. $\Lambda^{(2)}$).  
\item Every irreducible $\Uq$ module with highest weight belonging 
    to $\Lambda^{(1)}$ ( resp. $\Lambda^{(2)}$) is contained in 
    some repeated tensor products of ${\Bbb E}$ ( resp. ${\Bbb E}^\dagger$ )  
    as an irreducible component.   
\end{enumerate}  
\end{proposition}

More detailed structures of the covariant and contravariant 
tensor irreps can be understood, e.g.,  
their characters and super characters can be computed, 
the Clebsch - Gordan problem of irreps within a given tensor type 
can also be resolved by using the supersymmetric Young diagram method. 
Here we elucidate some general aspects of the Clebsch - Gordan 
problem, which will play an important role in the remainder 
of the paper.

Denote by $[\lambda]$ the equivalence class of irreps
with highest weight $\lambda$.  
For $\lambda$ and $\lambda'$ both belonging to $\Lambda^{(1)}$,  
we interpret  $[\lambda] + [\lambda']$ 
as the equivalence class of the direct sum representations, 
and $[\lambda]\cdot [\lambda']$ as that of the direct products. 
Let $\left[ \Lambda^{(1)} \right]$ 
be the $\Z^+$ module with a basis  
$\{ [\lambda] \ | \ \lambda\in \Lambda^{(1)} \}$.  
Then the `$\cdot$' operation  defines a multiplication on 
$\left[ \Lambda^{(1)} \right]$.
Clearly $[\lambda]\cdot [0]= [0]\cdot [\lambda]=[\lambda]$.    
Further more, from section V of \cite{Gould} we can 
deduce that if $[\lambda]\cdot [\lambda']  
=[\lambda^1] + [\lambda^2]+ ... + [\lambda^k]$, then none of 
the $\lambda^i$ is zero unless both $\lambda$ and 
$\lambda'$ are zero. 
This is in agreement with the fact that 
\ban \Lambda^{(1)}\bigcap \Lambda^{(2)} &=& \{0\}. \nan 
The dicussions above can be repeated word by word for the irreps 
with highest weights belonging to $\Lambda^{(2)}$.

\section{\normalsize QUANTUM SPECIAL LINEAR SUPERGROUP $GL_q(m|n)$}
For compact Lie groups in the classical setting, there exists the 
celebrated Tannaka - Krein duality theory \cite{Hewitt}, 
which enables the reconstruction of a group from the Hopf algebra of its 
representative functions.  The theory of quantum groups \cite{Faddeev} 
makes essential use of a quantum analog of the duality\cite{Woronowicz}, 
and is formulated entirely in terms of the algebra of functions. 
We will adopt the same philosophy here to formulate and study 
quantum supergroups.  However, we should mention that Lie 
supergroups are much more complicated than ordinary compact Lie 
groups in structures; at the best, the Tannaka - Krein duality holds 
in a restricted sense for Lie supergroups even at the classical 
situation, though we have not come across any treatment of the 
problem in the literature. 

\subsection{\normalsize Subalgebra of functions associated with the 
 vector irrep}
As before, we denote by $\pi$ the vector irrep of $\Uq$ relative to 
the standard basis $\{ v_a\ | \ a\in{\bf I}\}$ of  ${\Bbb E}$. Then  
\ban 
x v_a&=& \sum_{b} \pi ( x )_{b\, a} v_b, \ \ \ x\in \Uq. 
\nan 
Let $(\Uq)^0$ be the finite dual of $\Uq$. Consider the elements 
$t_{a\, b}$, $a, b\in{\bf I}$ of $(\Uq)^0$ satisfying 
\ban 
t_{a\, b} ( x ) &=& \pi ( x )_{a\, b}, \ \ \ \forall x\in \Uq.
\nan 
It is easy to show that the $t_{a\, b}$ indeed belong to $(\Uq)^0$.
Also note that $t_{a\, b}$ is even if $[a]+[b]\equiv 0 ( mod\, 2 )$, 
and odd otherwise.

Standard Hopf algebra theory asserts that $(\Uq)^0$ is a $\Zz$ 
- graded Hopf algebra with its structures dualizing those of $\Uq$.
Consider the subalgebra $G_q^{\pi}$ of $(\Uq)^0$ generated by 
$t_{a\, b}$, $a, b\in{\bf I}$. The multiplication which 
$G_q^{\pi}$ inherits from $(\Uq)^0$ is given by 
\ba 
\langle t\ t',\ x\rangle &=& \sum_{(x)} \langle t\otimes t', 
\ x_{(1)}\otimes x_{(2)} \rangle\nonumber\\ 
&=& \sum_{(x)} (-1)^{[t'][x_{(1)}]} \langle t,  x_{(1)} \rangle  
\langle t,\ x_{(2)} \rangle, 
\ \ \ \ \forall  t, t'\in G_q^{\pi}, \ x\in \Uq.  \label{product}   
\na 
To better understand the algebraic structure of $G_q^\pi$, 
we recall that the Drinfeld version of $\Uq$ admits a universal 
$R$ matrix, which in particular satisfies
\ban
R \Delta( x ) &=& \Delta'( x ) R,  \ \ \ \forall x\in \Uq.
\nan
Applying $\pi\otimes\pi$ to both sides of the equation yields 
\ba 
R^{\pi\, \pi} (\pi\otimes\pi)\Delta( x ) 
&=& (\pi\otimes\pi)\Delta'( x ) R^{\pi\, \pi}, \label{R}  
\na 
where 
\ban 
R^{\pi\, \pi} &:=&(\pi\otimes\pi)R\\ 
&=& q^{\sum_{a\in{\bf I}} e_{a\, a}\otimes e_{a\, a} (-1)^{[a]} } 
 + (q-q^{-1})\sum_{a<b} e_{a\, b}\otimes e_{b\, a} 
(-1)^{[b]} . 
\nan 
It is important to realize that equation (\ref{R})
 makes perfect sense within the Jimbo formulation of 
the quantized universal enveloping algebra $\Uq$, even when 
$q$ is specialized to a real parameter. 
We can re - interpret the equation in terms of the $t_{a\, b}$.
Set  $t=\sum_{a, b} e_{a\, b}\otimes t_{a\, b}$. Then
\ba
R^{\pi\, \pi}_{1 2}\, t_1\, t_2 
&=& t_2\, t_1\, R^{\pi\, \pi}_{1 2}.  \label{Faddeev}
\na   

The co - multiplication  $\Delta$ of $G_q^{\pi}$ is also defined 
in the standard way by 
\ban 
\langle \Delta( t_{a\, b} ), \ x\otimes y \rangle&=&
\langle t_{a\, b} , \ x y\rangle = \pi (x y )_{a\, b}, 
\ \ \ \ \forall x, y\in\Uq.  
\nan
We have 
\ba
\Delta( t_{a\, b} )&=&\sum_{c\in{\bf I}} 
(-1)^{([a]+[c])([c]+[b])} t_{a\, c}\otimes t_{c\, b}.   
\na 
$G_q^{\pi}$ also has the unit $\epsilon$, and the co - unit $1_{\Uq}$. 
Therefore, $G_q^{\pi}$ has the structures of a $\Zz$ - graded bi - 
algebra. However,  it does not admit an antipode, 
as will be explained later.

Let $\pi^{(\lambda)}$ be an arbitrary irreducible contravariant 
tensor representation of $\Uq$.  
We may also regard $\pi^{(\lambda)}$ 
as a representative of $[\lambda]$, where 
$\lambda$ $\in$ $\Lambda^{(1)}$.  
Define the elements $t^{(\lambda)}_{i\, j}$, 
$i, \, j= 1, 2, ..., dim_\C \pi^{(\lambda)}$, of $(\Uq)^0$ by 
\ban
t^{(\lambda)}_{i\, j}( x )&=& \pi^{(\lambda)}( x )_{i\, j}, 
\ \ \ \forall x\in\Uq. 
\nan 
It is an immediate consequence of Proposition \ref{tensor} that 
$t^{(\lambda)}_{i\, j}\in G_q^\pi$, for all $i,\ j$ and 
$\lambda\in\Lambda^{(1)}$, and every $f\in G_q^\pi$ can 
be expressed as a linear sum of these elements. 
{}From the representation theory of $\Uq$ we can deduce that 
these elements are also linearly independent. 
Introduce the vector spaces 
\ban 
T^{(\lambda)}=\bigoplus_{i, j=1}^{dim\pi^{(\lambda)}} 
\C t^{(\lambda)}_{i\, j}.
\nan 
Then  
\begin{proposition}  
As a vector space,  
\ban G_q^\pi&=&\bigoplus_{\lambda\in\Lambda^{(1)}} T^{(\lambda)}.\nan 
\end{proposition} 

Let us also denote the antipode of $(\Uq)^0$ by $S$. Then for any 
$t^{(\lambda)}_{i\, j}$ $\in G_q^\pi$,  
\ban 
S(t^{(\lambda)}_{i\, j} ) ( x ) &=& t^{(\lambda)}_{i\, j} ( S(x) ), 
\ \ \ \forall x\in\Uq. 
\nan 
That is, $S(t^{(\lambda)}_{i\, j} )$ are the matrix elements of the 
dual irrep of $\pi^{(\lambda)}$, the highest weight of which is 
not contained in $\Lambda^{(1)}$ unless $\lambda=0$.  Therefore, 
$G_q^\pi$ by itself does not admit an antipode.

\subsection{\normalsize Subalgebra of functions associated with the 
dual vector irrep}
Let $\{{\bar v}_a\ | \ a\in{\bf I}\}$ be the basis of ${\Bbb E}^\dagger$ 
dual to the standard basis of ${\Bbb E}$, i.e., 
\ban 
{\bar v}_a ( v_b )&=\delta_{a \, b}. 
\nan 
Denote by $\bar{\pi}$ the covariant vector irrep 
relative to this basis. Let ${\bar t}_{a\, b}$, $a, b\in{\bf I}$, 
be the elements of $(\Uq)^0$ such that 
\ban 
{\bar t}_{a\, b}( x )&=&{\bar\pi}( x )_{a\, b}, 
\ \ \ \forall x\in\Uq.
\nan  
Note that ${\bar t}_{a\, b}$ is even if $[a]+[b]\equiv 0 ( mod\,  2 )$, 
and odd otherwise.  These elements generate a $\Zz$ - graded bi - 
subalgebra $G_q^{\bar\pi}$ of $(\Uq)^0$ in the standard fashion. 
Here we merely point out that they obey the relation
\ba 
R^{{\bar\pi}\, {\bar\pi}}_{1 2}\ {\bar t}_1\ {\bar t}_2&=& 
{\bar t}_2\ {\bar t}_1 R^{{\bar\pi}\, {\bar\pi}}_{1 2}, 
\label{YB} \na  
where 
\ban 
{\bar t}&=&\sum_{a, b} e_{a\, b}\otimes {\bar t}_{b\, a}, \\ 
R^{{\bar\pi}\, {\bar\pi}}&=& ({\bar\pi}\otimes{\bar\pi})R,\\
&=& q^{\sum_{a\in{\bf I}} e_{a\, a}\otimes e_{a\, a} (-1)^{[a]} }
 + (q-q^{-1})\sum_{a>b} e_{a\, b}\otimes e_{b\, a} (-1)^{[b]}.
\nan
Also, the co - multiplication is given by 
\ban 
\Delta( \bar{t}_{a\, b} )&=&\sum_{c\in{\bf I}}
(-1)^{([a]+[c])([c]+[b])} \bar{t}_{a\, c}\otimes \bar{t}_{c\, b}.
\nan

Denote by ${\bar\pi}^{(-\bar{\lambda})}$ the irrep dual to 
${\pi}^{({\lambda})}$, $\lambda\in\Lambda^{(1)}$, in a given 
homogeneous basis. Introduce the elements 
$\bar{t}^{(-\bar{\lambda})}_{i\, j}$, 
$i, \, j=1, 2, ..., dim_{\C}{\pi}^{({\lambda})}$, 
of $(\Uq)^0$ such that 
\ban 
\bar{t}^{(-\bar{\lambda})}_{i\, j} ( x ) 
&=& {\bar\pi}^{(-\bar{\lambda})}_{i\, j} ( x ), \ \ \ \forall x\in\Uq. 
\nan 
Then it follows from Proposition \ref{tensor} that these elements 
form a basis of $G_q^{\bar\pi}$.  Set ${\bar T}^{(\mu)}$ 
$=\oplus_{i, j} \C \bar{t}^{(\mu)}_{i\, j}$, we have 
\begin{proposition}
\ban 
G_q^{\bar\pi}&=&\bigoplus_{\mu\in\Lambda^{(2)}} {\bar T}^{(\mu)}.
\nan 
\end{proposition}

\subsection{\normalsize Algebra $G_q$ of functions  on $GL_q(m|n)$}
We define the algebra $G_q$ of functions  on the quantum 
 general  linear supergroup $GL_q(m|n)$ to be the $\Zz$ - graded  
 subalgebra of $(\Uq)^0$ generated by $\{ t_{a\, b}, \ \bar{t}_{a\, b}$ 
$| \ a, b\in{\bf I}\}$.  The  $t_{a\, b}$ and 
$\bar{t}_{a\, b}$, besides obeying the the relations (\ref{Faddeev}) 
and (\ref{YB}), also  satisfy 
\ba R^{{\bar\pi}\, \pi}_{1 2}\ {\bar t}_1 \ t_2
&=& t_2 \ {\bar t}_1 \ R^{{\bar\pi}\, \pi}_{1 2}, \label{mix} \na
where 
\ban 
R^{{\bar\pi}\, \pi}  &:=& ({\bar\pi}\otimes \pi) R,\\ 
                     &=&
 q^{-\sum_{a\in{\bf I}} e_{a\, a}\otimes e_{a\, a} (-1)^{[a]} }
 - (q-q^{-1})\sum_{a<b} e_{b\, a}\otimes e_{b\, a}
(-1)^{[a]+ [b] + [a][b]}.
\nan 
Equation (\ref{mix}) enables us to factorize $G_q$ into 
\ba G_q&=&G_q^\pi \, G_q^{\bar\pi}.\label{factor}\na 
As both $G_q^\pi$ and $G_q^{\bar\pi}$ are $\Zz$ - graded bi - algebras, 
$G_q$ inherits a natural bi - algebra structure.  It also admits an 
antipode.  By considering  
\ban 
(x {\bar v}_a) (v_b)&=&(-1)^{[x][a]} {\bar v}_a ( S(x) v_b ), 
\ \ \ x\in \Uq,  
\nan 
where $\{ v_a\}$ is the standard basis of the vector irrep, and 
$\{{\bar v}_a \} $ is the basis of the covariant vector irrep 
dual to $\{ v_a\}$, we arrive at  
\begin{lemma} 
The antipode $S: G_q\rightarrow G_q$ is a linear anti - automorphism 
given by 
\ba 
S(t_{a\, b})&=& (-1)^{[a][b]+[a]} {\bar t}_{b\, a}, \nonumber\\
S({\bar t}_{a\, b})&=& (-1)^{[a][b]+[b]}  
q^{( 2\rho,\  \epsilon_a - \epsilon_b )} t_{b\, a} .
\na  
\end{lemma} 
Therefore, $G_q$ has the structures of a $\Zz$ - graded Hopf algebra.

Furthermore, $\ast$ - operations can also be constructed for $G_q$, 
thus turning it into a Hopf $\ast$ - superalgebra.  We have 
\ban 
\ast (t_{a\, b})&=& (-1)^{(\theta+[a])([a]+[b])} {\bar t}_{a\, b}, \\
\ast ({\bar t}_{a\, b})&=& (-1)^{(\theta+[a])([a]+[b])} t_{a\, b}. 
\nan 
where $\theta\in\Zz$.

An important property of $G_q$ is that it separates points of $\Uq$, 
that is, for any nonvanishing $x\in\Uq$,  there exists $f\in G_q$ 
such that $f( x )\ne 0$.  As a matter of fact, $G_q^\pi$ by itself
separates points of $\Uq$. Put differently, for any $u\in\Uq$, if 
$u\ne 0$, then $\pi^{\otimes p}(u)\ne 0$ for some $p\in\Z_+$.    

To verify our assertion, we first consider the corresponding 
proposition in the classical situation of $\U$ in detail.    
Let $E^{(0)}_{a\, b}$, $a, \, b\in{\bf I}$, be the 
standard generators of ${\frak g}$ embedded in its universal 
enveloping algebra.  In the vector irrep $\pi^{(0)}$, one 
has $$ \pi^{(0)}(E^{(0)}_{a\, b}) = e_{a\, b}.$$ 
We isolate the $u(1)$ subalgebra of ${\frak g}$ with the 
generator 
\ban 
Z^{(0)} &=& \sum_{a\in{\bf I}} E^{(0)}_{a\, a}, 
\nan 
and denote by $X_A^{(0)}$, $A=1, ...., (m+n)^2-1$, the elements
$E^{(0)}_{c\, c} - E^{(0)}_{c+1 \ c+1}$,  $c\in{\bf I}'$, and   
$E^{(0)}_{a\, b}$, $a\ne b$,  in any fixed ordering. Then a PBW basis 
for $\U$ is given by, 
$$ \{ B^{(0)}_{k,\ A_1 ... A_l}= (Z^{(0)})^k X^{(0)}_{A_1} 
... X^{(0)}_{A_l}\, |\, 
k, l\in\Z_+, \ A_i\le A_{i+1}, \ A_i\ne A_{i+1} 
\ \mbox{if}\  [X^{(0)}_{A_i}]=1\}.$$ 
Set $\pi^{(0)}(X^{(0)}_{A})=e_A$. Denote by ${\cal M}$ the vector space 
of $(m+n)\times (m+n)$ matrices, and define 
\ban 
{\cal R}_k&=&\sum_{i=0}^{k-1} \underbrace{{\cal M}\otimes ... 
\otimes {\cal M}}_i  \otimes I \otimes 
\underbrace{{\cal M}\otimes ... \otimes {\cal M}}_{k-1-i}. 
\nan  
Let 
\ban 
b_{A_1 \, ...\, A_k}&=&\sum_{\sigma\in S_k} (-1)^{|\sigma_{\{A\}}|} 
e_{A_{\sigma(1)}}\otimes e_{A_{\sigma(2)}}\otimes ... \otimes 
e_{A_{\sigma(k)}},  
\nan   
where $|\sigma_{\{A\}}|$ is the number of permutations required 
 amongst odd elements in order to change $X_{A_1}\otimes X_{A_2} 
\otimes ... \otimes X_{A_k}$ to $X_{A_{\sigma(1)}} 
\otimes X_{A_{\sigma(1)}}\otimes ... \otimes X_{A_{\sigma(k)}}$. 
Clearly, the elements 
$$\{ b_{A_1 \, ...\, A_k}\, |\, k\in\Z_+, 
\ A_i\le A_{i+1}, \ A_i\ne A_{i+1} 
\ \mbox{if}\  [X_{A_i}]=1\}$$ 
are linearly independent in ${\cal M}^{\otimes k}$, and we will 
denote by ${\cal L}_k$ their linear span.  By considering the 
trace ( not the supertrace! ) on each factor of ${\cal M}^{\otimes k}$, 
we can easily see that ${\cal L}_k$ 
intersects ${\cal R}_k$ trivially. 
Therefore, 
\ban 
(\pi^{(0)})^{\otimes (k+p)}( B^{(0)}_{0, \ A_1 ... A_k} ) &=& 
\left( (\pi^{(0)})^{\otimes k}\otimes (\pi^{(0)})^{\otimes p}\right) 
( B^{(0)}_{0, \ A_1 ... A_k} ) \\
&=& b_{A_1 \, ...\, A_k}\otimes I^{\otimes p} + r_{k, \, p}, 
\ \ \ r_{k, \, p}\in {\cal R}_k\otimes {\cal M}^{\otimes p}, 
\nan 
are linearly indepenedent as elements of ${\cal M}^{\otimes (k+p)}$.

Consider  $u\in\U$ given by 
\ban 
u&=&\sum_{k=0}^K \sum_{l=0}^L \sum_{\{A\}} 
C_{k,\ A_1 ... A_l} B^{(0)}_{k,\ A_1 ... A_l}, 
\ \ \ \ C_{k,\ A_1 ... A_l}\in\C. 
\nan 
Using 
$$(\pi^{(0)})^{\otimes p}(Z^k)=p^k\ I^{\otimes p},$$ 
we immediately see that $(\pi^{(0)})^{\otimes p}(u)=0$, $\forall 
p> L$, requires 
\ban 
\sum_{k=0}^K p^k C_{k,\ A_1 ... A_l} &=&0, \ \ \ 
\forall p> L,    
\nan 
which forces all the $C_{k,\ A_1 ... A_l}$ to vanish. This 
completes the proof for the classical case.\\ 

\noindent
{\bf Remark}: There is something slightly unnatural about our 
proof, that is, the combination $E^{(0)}_{m\, m} - E^{(0)}_{m+1\ m+1}$ does 
not belong to $sl(m|n)\subset \gl$, and this in turn forced us  
to consider the ordinary trace instead of the supertrace in proving 
${\cal L}_k\cap {\cal R}_k $ $=\{ 0 \}$.  We can avoid this unnaturalness 
when $m\ne n$ by using $E^{(0)}_{m\, m}+E^{(0)}_{m+1\ m+1}$ instead, 
but not when $m=n$.\\

With the above preparations we can now readily prove our assertion 
for the quantum superalgebra.  We first consider the Drinfeld 
version of $\Uq$.  Similar to the  classical case, we set 
\ban
Z&=& \sum_{a\in{\bf I}} E_{a\, a},
\nan
and denote by $X_A$, $A=1, ...., (m+n)^2-1$, the elements
$E_{c\, c} - E_{c+1 \ c+1}$, $c\in{\bf I}'$ and  
$E_{a\, b}$, $a\ne b$,  in a fixed ordering. Then  
$$ \{ B_{k,\ A_1 ... A_l}= Z^k X_{A_1}
... X_{A_l}\, |\,
k, l\in\Z_+, \ A_i\le A_{i+1}, \ A_i\ne A_{i+1}
\ \mbox{if}\  [X_{A_i}]=1\}$$
forms a PBW basis for $\Uq$. 
Given  
\ban 
u=\hbar^k ( u_0 + \hbar u_1 + \hbar^2 u_2 + ... ).  
\nan 
Each $u_i$ is  a finite $\C$ - combination of 
some $ B_{k,\ A_1 ... A_l}$, and $u_0$ is assumed to be 
nonzero.   Then it follows from the classical case that 
there exist infinitely many $p\in\Z_+$ such that 
\ban 
\pi^{\otimes p} ( u) \not\equiv 0  ( mod\, \hbar^{k+1} ). 
\nan 

For the Jimbo algebra, we observe that ordered monomials in 
$E_{a\, b}$, $a\ne b$, and $K_a^{\pm 1}$ form a basis of $\Uq$.  
Given $u\in\Uq$, and a positive integer $p$, we consider 
the matrix elements of $\pi^{\otimes p}(u)|_{q=\exp(\hbar)}$ as 
power series in $\hbar$.  $\pi^{\otimes p}(u)\ne 0$ if and only 
if some of these power series do not vanish identically. 
Now for the purpose of computing $\pi^{\otimes p}(u)|_{q=\exp(\hbar)}$, 
we can make the identification 
\ban 
\pi^{\otimes p}(K_a)&=&\sum_{k=0}^\infty 
{{(-1)^{k [a]} \hbar^k }\over{k!} } e_{a\, a}(p)^k,\\  
e_{a\, a}(p)&=& \sum_{i=0}^{p-1} \underbrace{I\otimes ... \otimes  I}_i
 \otimes e_{a\, a}\otimes \underbrace{I\otimes ...\otimes I}_{p-i-1}. 
\nan 
This takes us back to the Drinfeld algebra situation, and we 
have already shown that in that situation the $\pi^{\otimes p}$, 
$p\in\Z_+$, separates points of $\Uq$.

We summarize the discussions of this section into a Proposition, 
points ii) and iii) of which may be considered as 
a partial generalization of the classical
Peter - Weyl theorem to the quantum supergroup
in an algebraic setting: 
\begin{proposition} i) $G_q$ is a  $\ast$ - Hopf superalgebra; \\ 
ii).  $G_q$ separates points of $\Uq$;\\ 
iii). The following elements span $G_q$:
\ban
 t^{(\lambda)}_{i\, j} \ {\bar t}^{(\mu)}_{i'\, j'},  & \ &
i,\ j = 1, 2, ..., dim\pi^{(\lambda)}, \ \
 \lambda\in\Lambda^{(1)},\\
& \ & i',\ j' = 1, 2, ..., dim{\bar\pi}^{(\mu)}, \ \
 \mu\in\Lambda^{(2)}. 
\nan
\end{proposition}
However, we should point out that these elements are {\em not} linearly
independent.

\section{\normalsize INDUCED REPRESENTATIONS OF $G_q$  }
We will develop parabolic induction for representations of 
$G_q$ in this section.  Recall that corresponding to every locally 
finite right co - module $\omega: W$  $\rightarrow W\otimes G_q$ 
over $G_q$, there exists a unique left $\Uq$ module 
$\Uq\otimes W$ $\rightarrow W$  with the module action defined by 
\ban 
x\, w &=& \omega(w)(x), \ \ \ x\in\Uq, \ \ w\in W. 
\nan 
A similar correspondence exists for left $G_q$ co - modules and 
right $U_q$ modules.  
Therefore,  we can describe the representation theory of $G_q$ 
in both the $G_q$ co - module language and $\Uq$  module language,
depending on which one is more convenient in a given situation. 
We will largely use the latter here.

\subsection{\normalsize Parabolic subalgebras of $U_q(gl(m|n))$} 
Let $\bf\Theta$ be a subset of ${\bf I^\prime}$.  
Introduce the following sets of elements of $\Uq$:
\ban
{\cal S}_l&=&\{ K_a^{\pm 1}, a\in {\bf I}; \
\ E_{c\ c+1}, \ E_{c+1\ c},  \ c\in {\bf\Theta}\};\\
{\cal S}_{p_+}&=& {\cal S}_l \cup \{ E_{c\ c+1},  
     \ c\in {\bf I^\prime} \backslash {\bf\Theta}\}; \\ 
{\cal S}_{p_-}&=& {\cal S}_l \cup \{ E_{c+1\ c},  
     \ c\in {\bf I^\prime} \backslash {\bf\Theta}\}.
\nan
The elements of each set generate a $\Zz$  - graded Hopf 
subalgebra of $\Uq$.   We denote by $\Ul$ the Hopf subalgebra 
generated by the elements of ${\cal S}_l$, 
and by $U_q({\frak p}_\pm)$ the Hopf subalgebras 
respectively generated by the elements of ${\cal S}_{p_\pm}$.
In the classical limit, the Hopf subalgebras $U_q({\frak p}_\pm)$ 
coincide with the universal enveloping
algebras of  parabolic subalgebras  of the
Lie superalgebra $\gl$.  Therefore, we will call 
$U_q({\frak p}_\pm)$ parabolic subalgebras of $\Uq$.   

Let $V_\mu$ be a finite dimensional irreducible $\Ul$ module. 
Then $V_\mu$ is of highest weight type.  Let $\mu$ be the  
highest weight and $\tilde\mu$ the lowest weight of $V_\mu$ 
respectively.   
We can extend $V_\mu$ in a unique fashion to a $\Up$ module, 
which we still denote by  $V_\mu$, 
such that the elements of ${\cal S}_{p_+}\backslash {\cal S}_l$ 
act by zero.  Similarly, $V_\mu$ also leads to a $\Um$ 
module, on which the elements of ${\cal S}_{p_-}\backslash {\cal S}_l$  
act by zero.  It is not difficult to see that all 
finite dimensional irreducible $\Upm$ modules are of this kind. 
 
Consider a finite dimensional irreducible $\Uq$ module 
$W(\lambda)$ with highest weight $\lambda$  and lowest weight 
$\bar\lambda$.  $W(\lambda)$ can be restricted into a 
$\Up$ or $\Um$ module in a natural way, and the resultant 
module is always indecomposable, but not irreducible in general.     

Consider first the case of $\Up$.  We wish to examine 
the $\Zz$ - graded vector space $Hom_{\Up} ( W(\lambda), \ V_\mu )$, 
which graded - commutes with $\Up$, namely, 
\ban 
p\, \phi - (-1)^{[p][\phi]} \phi\, p &=&0, \ \ \ p\in\Up, 
\ \phi\in Hom_{\Up} ( W(\lambda), \ V_\mu ).  
\nan 
Because of the 
irreducibility of $V_\mu$, every non - zero $\phi$ 
$\in Hom_{\Up} ( W(\lambda), \ V_\mu )$ must be surjective, 
and thus $V_\mu\cong W(\lambda)/ Ker\phi$.    
As a $\Up$ module, $W(\lambda)$ is indecomposable, 
and contains a unique maximal proper submodule $M$ such 
that the lowest weight vector $w_-$ of $W(\lambda)$ 
does not belong to $M$.  Therefore, $Ker\phi = M$, 
and $V_\mu=\phi ( \Ul w_- )$. 
This forces  $\bar\lambda = \tilde\mu$, and all 
  elements of $Hom_{\Up} ( W(\lambda), \ V_\mu )$ 
are scalar multiples of one another.  It is worth observing that 
the map $\phi$ may be odd. In fact its degree is given by 
$[\phi]\equiv [w_-]+[\phi(w_-)]\, ( mod 2 )$.  
The case of $\Um$ can be studied in exactly the same way. 
To summarize, we have 
\begin{lemma}\label{Hom} 
\ban
\dim_{\C} \mbox{Hom}_{\Up} ( W(\lambda), \ V_\mu) &=& 
    \left\{ \begin{array}{l l}
           1, &  \bar\lambda={\tilde\mu}, \\
           0, &  \bar\lambda\ne {\tilde\mu}.
     \end{array}\right.\\
\dim_{\C} \mbox{Hom}_{\Um} ( W(\lambda), \ V_\mu) &=& 
    \left\{ \begin{array}{l l}
           1, &  \lambda=\mu, \\
           0, &  \lambda\ne \mu.
     \end{array}\right.
\nan
\end{lemma}

\subsection{\normalsize Induced representations and quantum  superbundles} 
Let us first introduce  two types of left actions of $\Uq$ on $G_q$, which
correspond to the left and right translations in the classical
situation.   

Define a bilinear map $\cdot:  \Uq \otimes G_q \rightarrow G_q$ by
\ba
x\otimes f &\mapsto&  x\cdot f\nonumber \\
    &=&\sum_{(f)} \langle f_{(1)}, \ S^{-1}(x) \rangle f_{(2)}, 
\na
which can be easily shown to satisfy
\ban
(x\cdot f)(y)&=& (-1)^{[x][y]} f( S^{-1}(x) y ),\\ 
x\cdot(y\cdot f)&=& (x y)\cdot f, \ \ \ \ x,\ y\in \Uq, \ f\in G_q.  
\nan
(We assume that the elements $x,\ y\in \Uq$  and $g,\ f\in
G_q$ are homogeneous for the sake of simplicity.  All the statements
below  generalize to inhomogeneous elements in the obvious way.)
Therefore, this defines a left action of $\Uq$ on $G_q$, 
which corresponds to the left translation of Lie groups in the
classical situation. It is worth observing that we may replace
$S^{-1}$ in the above definition, and arrive at a different  
left action.

Another left action `$\circ$' of $\Uq$ on $G_q$ can be defined by 
\ba 
x\circ f&=&\sum_{(f)} f_{(1)} (-1)^{[x]([f]+[x])}\ \langle f_{(2)}, 
\ x\rangle. \label{circ}  
\na 
Straightforward calculations can  show that 
\ban y\circ ( x\circ f ) &=& (x y)\circ f;\\ 
(x\circ f) ( y ) &=& f(y x), \\ 
(id_{G_q}\otimes x\circ) \Delta(f) &=& \Delta( x\circ f ).  \nan 
This corresponds to the right translation in the classical theory. 
It graded - commutes with the action `$\cdot$',  namely, 
\ban x\circ( y\cdot f) &=& (-1)^{[x][y]} y\cdot(x\circ f). 
\nan 

Let $\Uqp$ denote either $\Up$ or $\Um$.  Given any  
finite dimensional left $\Uqp$ module $V$,  
we form the tensor product $ V \otimes_{\C} G_q$, which is 
a subspace of functions $\Uq\rightarrow V$:   
\ban 
\zeta &=&\sum v_i \otimes f_i  \in V\otimes G_q ,\\ 
 x&\in& \Uq,  \\  
\zeta (x)&=&\sum f_i (x) v_i. 
\nan 
The left actions `$\cdot$' and `$\circ$' of $\Uq$ on $G_q$ 
can be extended in an obvious way to actions on $V\otimes_{\C} G_q$
\ban 
x\cdot\zeta &=& \sum (-1)^{[x][v_i]} v_i \otimes x\cdot f_i, \\
x\circ\zeta&=& \sum (-1)^{[x][v_i]} v_i \otimes x\circ f_i,  
\ \ \ \  x\in\Uq. 
\nan   
Furthermore, there also exists a co - action  $\omega$ 
of $G_q$ on $V\otimes_{\C} G_q$ defined by 
$\omega=id_V\otimes \Delta'$, where $\Delta'$ represents the 
opposite co - multiplication of $G_q$. 

Consider the subspace of $V\otimes_{\C} G_q$ defined by   
\ban
{\cal O}^V &=& \{ \zeta\in V\otimes_{\C} G_q \, |\, 
p\circ\zeta  = (\, S(p) \otimes id_{G_q}) \zeta,  \ \forall 
\ p\in\Uqp\}.   
\nan 
\begin{lemma} 
${\cal O}^V$ furnishes a left $\Uq$ module under `$\cdot$', 
and at the same time a right $G_q$ co - module under $\omega$. 
\end{lemma} 
{\em Proof}: The Lemma can be confirmed
by direct calculations.  For $x\in\Uq$, $p\in\Uqp$, $\zeta\in{\cal O}^V$,
we have
\ban
p\circ  ( x\cdot\zeta )&=& (-1)^{[x][p]}x\cdot(p\circ\zeta)\\   
&=&(\, S(p) \otimes  id_{G_q}) ( x\cdot\zeta ); \\
(\, p\circ\otimes  id_{G_q}) \omega(\zeta) 
&=& (\, p\circ\otimes  id_{G_q}) ( \, id_V \otimes \Delta' )\zeta\\
&=&(\, id_V \otimes \tau ) ( \, id_V\otimes id_{G_q}\otimes p\circ) 
( \, id_V \otimes \Delta) \zeta\\ 
&=&(\, id_V \otimes \tau ) ( \, id_V \otimes \Delta' )(p\circ \zeta)\\  
&=&\omega(\, S(p) \otimes id_{G_q})\zeta, 
\nan  
where $\tau$ is the flip mapping.  

We call ${\cal O}^V$ the induced $\Uq$ module, and  also  
the induced $G_q$ co - module, which gives rise to 
a representation of $G_q$.   
A conceptual understanding of ${\cal O}^V$ can be gained by 
considering its classical analog.  Let $P$ be a parabolic 
subgroup of the complex Lie supergroup $SL(m|n)$, and $E$ 
a finite dimensional representation of $P$. Then 
$SL(m|n)\times_P E$, the quotient space of $SL(m|n)\times E$ 
under the equivalence relation $(g, v) \sim ( g p,  p^{-1} v)$ 
for all $p\in P$, defines a super vector bundle over the 
supermanifold $SL(m|n)/P$.  A function $f: SL(m|n)\rightarrow 
E$ satisfying $f(g p) = p^{-1} f(g)$, $\forall p\in P$ 
defines a section of the bundle $s_f: SL(m|n)/P \rightarrow$ 
$SL(m|n)\times_P E$.  
Analogously, we may regard ${\cal O}^V$ as the vector space 
of sections of a quantum super vector bundle over the 
quantum counter part of $SL(m|n)/P$.   

It is of great importance to systematically develop the theory 
of quantum homogeneous super vector bundles, and the subject 
will be investigated to depth in a forthcoming publication. 
In this paper, we will restrict ourselves to issues directly 
related to representation theory, and will not further ponder on   
noncommutative geometry, except for the last section, where we 
will discuss in some detail quantum projective superspaces 
when dealing with explicit realizations of the 
skew supersymmetric tensor irreps and their duals.

We have the following 
quantum analog of Frobenius reciprocity.  
\begin{proposition} Let $W$ be a quotient $\Uq$ module of 
$\bigoplus _{k, l=0}^\infty$ 
${\Bbb E}^{\otimes k}\otimes ({\Bbb E}^*)^{\otimes l}$ ( the restriction 
of which furnishes a $\Uqp$ module in a natural way.).  Then  
there is a canonical isomorphism 
\ba 
\mbox{Hom}_{\Uq} ( W, \ {\cal O}^V ) 
&\cong& \mbox{Hom}_{\Uqp} ( W, \ V), 
\na
\end{proposition} 
{\em Proof}:  We prove the Proposition by explicitly constructing 
the isomorphism, which we  claim to be  the linear map 
\ban 
F: \mbox{Hom}_{\Uq} ( W, \ {\cal O}^V ) &\rightarrow& 
   \mbox{Hom}_{\Uqp} ( W, \ V), \\ 
    \psi &\mapsto& \psi(1_{\Uq}), 
\nan 
with the inverse map 
\ban 
\bar{F}: \mbox{Hom}_{\Uqp} ( W, \ V) &\rightarrow& 
        \mbox{Hom}_{\Uq} ( W, \ {\cal O}^V ), \\ 
       \phi &\mapsto& \bar{\phi}, 
\nan 
where $\bar{\phi}$ is defined by
\ban 
\bar{\phi}(w)(x)&=&(-1)^{[x]([w]+1)} \phi(S(x) w),  \ \ \ 
 x\in\Uq, \ w\in W.
\nan  

As for $F$, we need to show that its image is contained in 
$\mbox{Hom}_{\Uqp} ( W, \ V)$. This is indeed the case, as 
\ban 
p ( F \psi (w) ) &=&(p\cdot\psi(w)) (1_{\Uq})\\
                      &=&(-1)^{[\psi][p]}  F \psi (p w), \ \ \ \  
p\in\Uqp, \ \ w\in W.  
\nan

In order to show that ${\bar F}$ is the inverse of 
$F$, we first need to demonstrate that the image $Im({\bar F})$ 
of ${\bar F}$ is contained in  
$\mbox{Hom}_{\Uq} ( W, \ {\cal O}^V )$. 
Note that $Im({\bar F})$ $\subset$  $\mbox{Hom}_{\C}  
( W, \ V \otimes G_q)$,  since
$W$ is a subquotient of $\bigoplus _{k, l=0}^\infty$
${\Bbb E}^{\otimes k}\otimes ({\Bbb E}^*)^{\otimes l}$. 
Some relatively simple manipulations  lead to  
\ban 
(y\cdot{\bar\phi} (w) )(x)  &=&(-1)^{[y][{\bar\phi}]+[x]([w]+[x]+[y])}
                      \phi( S(x) y w )\\
                      &=& (-1)^{[y][{\bar\phi}]} {\bar\phi} ( y w ) (x),\\ 
( p\circ {\bar\phi} (w) ) (x) 
               &=&(-1)^{[x]([w]+1)+[p][\phi]} \phi( S(p)  S(x) w)\\ 
               &=& S(p) ( {\bar\phi} (w) (x) ), \ \ \ \  
x, \  y\in\Uq, \ \ p\in\Uqp, \ \ w\in W. 
\nan 
Therefore,  $Im({\bar F})$ $\subset$ 
$\mbox{Hom}_{\Uq} ( W, \ {\cal O}^V )$.

Now we show that $F$ and ${\bar F}$ are inverse to each other. 
For $\psi\in$ $\mbox{Hom}_{\Uq} ( W, \ {\cal O}^V )$, and 
$\phi$ $\in$ $\mbox{Hom}_{\Uqp} ( W, \ V)$, we have  
\ban 
(F \bar{F}  \phi)(w) &=& (\bar{F} \phi)(w) (1_{\Uq})\\
                             &=& \phi(w), \\
(\bar{F} F \psi)(w)(x)&=& (-1)^{[x]([w]+1)} ( F \psi) (S(x) w)\\ 
        &=&(-1)^{[x]([w]+1)} \psi (S(x) w) (1_{\Uq})\\
        &=&(-1)^{[x]([\psi(w)]+1)}  (S(x)\cdot \psi (w)) (1_{\Uq})\\
        &=& \psi (w) (x), \ \ \ \ x\in\Uq, \ \ w\in W. 
\nan 
This completes the proof of the Proposition.

\subsection{\normalsize Quantum Borel - Weil theorem 
for covariant and contravariant tensor irreps} 
In this subsection we study in detail the covariant and contravariant 
tensor irreps of $\Uq$ within the framework of parabolic induction. 
Our main result here will be a quantum version of the 
Borel - Weil theorem for these irreps. 

Let $V$ be a finite dimensional 
irreducible $\Uqp$ module, with the $\Ul$  highest  weight 
$\mu$ and $\Ul$ lowest  weight $\tilde\mu$. 
For the purpose of studying the 
tensor representations, we only need to consider 
\ban
{\overline{\cal O}} (\mu) &=& {\cal O}^V \cap
                          \left( V\otimes G_q^\pi\right), \\
{\cal O} (\mu) &=& {\cal O}^V \cap
                         \left(V\otimes G_q^{\bar\pi}\right).
\nan
 
Let us study ${\cal O} (\mu)$ first.  
A typical element of ${\cal O} (\mu)$ is of the form
\ban
\zeta &=& \sum_{\lambda\in\Lambda^{(1)}}
          \sum_{\alpha, \beta, i} c_{\alpha\, \beta, \ i}^\lambda
          \ v_i\otimes{\bar t}_{\alpha\, \beta}^{(\lambda^\dagger)},
\nan
where $\{ v_i \}$ is a basis of $V$, and the
$c_{\alpha\, \beta, \ i}^\lambda$ are complex numbers. 
The ${\bar t}_{\alpha\, \beta}^{(\lambda^\dagger)}$ 
are elements of the Peter - Weyl basis for $G_q^{\bar\pi}$, which, 
needless to say, are polynomials in ${\bar t}_{a b}$, 
$a, \, b\in{\bf I}$. 
The property that $(p\circ\zeta)$ $=$  $(\, S(p)\otimes id_{G_q}) \zeta$, 
$\forall p\in\Uqp$  leads to
\ba
\sum_{\gamma, i} (-1)^{[p]([\gamma]+[v_i])}
     c_{\alpha\, \gamma, \  i}^\lambda\ 
{t}_{\gamma\, \beta}^{(\lambda)}(p) v_i
&=& \sum_{i} c_{\alpha\, \beta, \ i}^\lambda\  p\, v_i,
\ \ \ \ \forall p\in\Uqp.\label{morphism} 
\na

Let $W(\lambda)$ with the basis  $\{w_\alpha\}$ be the irreducible 
$\Uq$ module associated with the irrep $t^{(\lambda)}$. 
We define the linear maps between $\Zz$ graded {\em vector spaces }  
\ban 
\phi_\lambda^{(\alpha)}:\  W(\lambda)& \rightarrow& V,\\ 
 w_\beta& \mapsto& \sum_i c_{\alpha\, \beta, \ i}^\lambda \, v_i.     
\nan 
There is no particular significance attached to the maps at this stage, 
apart from the mere fact that they can be employed to re - express  
equation (\ref{morphism}) as
\ban 
\sum_\gamma (-1)^{[p][\phi_\lambda^{(\alpha)}]}  
  {t}_{\gamma\, \beta}^{(\lambda)}(p) 
  \phi_\lambda^{(\alpha)} (w_\gamma)  
&=& p\, \phi_\lambda^{(\alpha)} (w_\beta). 
\nan 
We emphasize that this equation is entirely equivalent to 
(\ref{morphism}).  Now something of crucial importance appears: 
this equation requires that each $\phi_\lambda^{(\alpha)}$ 
be a $\Uqp$ module homomorphism of degree $[\phi_\lambda^{(\alpha)}]$.   
Lemma \ref{Hom}  forces 
\ban  \phi_\lambda^{(\alpha)} = c_\alpha \, \phi_\lambda, &  
c_\alpha\in \C, \nan   
and  $\phi_\lambda$ may be nonzero only when 
\ban 
i). &    \bar\lambda={\tilde\mu},   & \mbox{if}\  \Uqp=\Up, \\
ii). &    \lambda=\mu, & \mbox{if}\  \Uqp=\Um. 
\nan 
In these cases, ${\cal O} (\mu)$ is spanned by  
\ban 
\zeta_\alpha &=& \sum_{\beta} \phi_\lambda(w_\beta)\otimes 
 {\bar t}_{\alpha\, \beta}^{(\lambda^\dagger)},    
\nan 
which are  obviously linearly independent. Furthermore, 
\ba 
x\cdot\zeta_\alpha &=& (-1)^{[x][\phi_\lambda]} 
 \sum_{\beta} t^{(\lambda)}_{\beta\, \alpha} (x)\  \zeta_{\beta}, 
\ \ \ \ \ x\in\Uq.  \label{irrep} \na   
The case of ${\overline{\cal O}}(\mu)$ can be studied in exactly the
same way.  To summarize, we have the following quantum analog of 
Borel - Weil theorem for the covariant and contravariant tensor irreps  
\begin{proposition} 
As  $\Uq$ modules, 
\ba 
{\cal O} (\mu)&\cong&
 \left\{ \begin{array}{l l l } 
  W((-{\tilde\mu})^\dagger), & \mbox{if} \ {\tilde\mu}\in -\Lambda^{(2)},  
           & \Uqp=\Up, \\
  W(\mu), & \mbox{if} \ \mu\in \Lambda^{(1)}, 
           & \Uqp=\Um,\\  
   \{0\},  &\mbox{otherwise}. 
                     \end{array}\right. \\ 
{\overline{\cal O}}(\mu)&\cong&
 \left\{ \begin{array}{l l l }
  W((-{\tilde\mu})^\dagger), & \mbox{if} \ {\tilde\mu}\in -\Lambda^{(1)},
           & \Uqp=\Up, \\
  W(\mu), & \mbox{if} \ \mu\in \Lambda^{(2)},
           & \Uqp=\Um,\\
   \{0\},  &\mbox{otherwise}.
                     \end{array}\right.
\na 
\end{proposition} 
In the Proposition, the notation $W(\lambda)$  
signifies the irreducible 
$\Uq$ module with highest weight $\lambda$.

\section{\normalsize QUANTUM PROJECTIVE SUPERSPACES, \ \ \    
SKEW SUPERSYMMETRIC TENSOR IRREPS AND THEIR DUALS}
We will apply the general theory developed in the last section 
to study two infinite classes of irreps, namely, the 
skew supersymmetric tensor irreps and their dual irreps. 
Explicit realizations of these irreps will be given 
in terms of sections of quantum super vector bundles over 
quantum projective superspaces. 
\subsection{\normalsize Quantum projective superspaces}  
Let $U_q({\frak g}')$, ${\frak g}'=gl(m|n-1)$,
 be the subalgebra of $\Uq$ generated by the following elements  
\ban 
\{ K_a,  \ a\in{\bf I}';  E_{c\, c+1}, 
\ E_{c+1\ c}, \ \ c\in{\bf I}'\backslash \{m+n-1\} \}.  
\nan 
Clearly $U_q({\frak g}')$ is a Hopf subalgebra. 
Define 
\ban 
{\cal A}_+&=&\{ f\in G_q^\pi \, |\ f(x p) = \epsilon ( p ) f ( x ), \ 
             \forall x\in\Uq,\ p\in U_q({\frak g}')\}, \\  
{\cal A}_- &=& \{ f\in G_q^\pi \, |\ f(x p) = \epsilon ( p ) f ( x ), 
             \ \forall x\in\Uq,\ p\in U_q({\frak g}')\}.    
\nan 
The Hopf algebra structure of $U_q({\frak g}')$ implies that both 
${\cal A}_+$  and ${\cal A}_- $ are subalgebras of $G_q$. 
Together they  generate another subalgebra of $G_q$, 
which we will denote by $\Sq$. 
Set
\ban 
z_a = t_{a\ m+n}, &  {\bar z}_a = {\bar t}_{a\ m+n}, & a\in{\bf I}.  
\nan 
Then $z_a$ and ${\bar z}_a$ are conjugate to each other under the 
$\ast$ - operation with $\theta=0$. More explicitly, 
  $$*(z_a)={\bar z}_a, \ \ \ \ \forall a\in{\bf I}.$$  
Now $\Sq$ is generated by the $z$'s and ${\bar z}$'s,   
which satisfy the following commutations relations  
\ban 
z_a\,  z_b &=& (-1)^{[z_a] [z_b]}\,  q\,  z_b\, z_a, \ \ \ a<b, \\
(z_c)^2 &=& 0, \ \ \ c\le m; \\
{\bar z}_a {\bar z}_b &=& (-1)^{[{\bar z}_a][{\bar z}_b]}
q^{-1} \, {\bar z}_b {\bar z}_a, \ \ \ a<b, \\
({\bar z}_c)^2&=& 0, \ \ \ c\le m;\\  
{\bar z}_a\, z_b &=& q (-1)^{[{\bar z}_a][z_b]} z_b \, {\bar z}_a 
+ \delta_{a\, b}\left\{ (1-q_a^{-1}){\bar z}_a\, z_a\right.\\  
    &-&\left. (-1)^{[{\bar z}_a]} (q-q^{-1})\sum_{c<a} {\bar z}_c\, z_c 
    \right\}, \ \ \ \forall a,\, b\in{\bf I}, \\
\sum_{c\in{\bf I}} {\bar z}_c\, z_c &=& 1.  
\nan 
It can be shown that the last two equations imply that 
\ban 
\sum_{c\in{\bf I}} q^{(2\rho, \ \epsilon_c)} \, z_c \, {\bar z}_c
       &=& q^{(2\rho, \ \epsilon_{m+n})}. 
\nan 
$\Sq$ furnishes a right $G_q$ co - module algebra,  with the co - module 
action $\omega:  \Sq$ $ \rightarrow \Sq \otimes G_q$ defined by 
\ban 
\omega( z_a )&=&\sum_{c\in{\bf I}}
z_c \otimes  t_{a\, c}, \\ 
\omega( \bar{z}_a )&=&\sum_{c\in{\bf I}}
\bar{z}_c \otimes \bar{t}_{a\, c}.
\nan 
Also, $\Sq$ gives rise to a right $\Uq$ module algebra 
with the module action `$\circ$' defined by (\ref{circ}).  This module 
algebra structure restricts naturally to a module algebra structure 
over $U_q({\frak g}')\otimes U_q(gl(1))$,  where  $U_q(gl(1))$ 
is generated by $K_{m+n}^{\pm 1}$.   The action of 
$U_q({\frak g}')$ on $\Sq$ is trivial as following the 
definitions of ${\cal A}_\pm$;  $U_q(gl(1))$ also acts in a 
very simple manner. To be explicit, we introduce the notations that 
for   $L=(\theta_1, ..., \theta_m; l_1, ..., l_n)$ $\in $ 
$\{0, \ 1\}^{\otimes m} \otimes \Z^{\otimes n}_+$,       
$|L|=\sum_{i=1}^m \theta_i + \sum_{\mu=1}^n l_\mu$. 
Set 
\ba 
Z^L &=& z_1^{\theta_1} ... z_m^{\theta_m} z_{m+1}^{l_1} ... z_{m+n}^{l_n}, 
\nonumber\\ 
{\overline Z}^L &=& {\bar z}_1^{\theta_1} ... {\bar z}_m^{\theta_m}
       {\bar z}_{m+1}^{l_1} ... {\bar z}_{m+n}^{l_n}. 
\label{Z}
\na 
Then for any $k\in\Z$, and $p\in U_q({\frak g}')$, we have 
\ba 
(p K_{m+n}^k )\circ ( Z^L\, {\overline Z}^{L'} ) 
&=& \epsilon ( p ) q^{k(|L'| - |L|)} Z^L\, {\overline Z}^{L'}. 
\label{U}
\na 

We will define the quantum projective superspace $\CP$ to be the 
$U_q(gl(1))$ invariant subalgebra of $\Sq$, namely, 
\ba 
\CP&=& \left( \Sq\right)^{U_q(gl(1))}. 
\na

\subsection{\normalsize  Skew supersymmetric 
            tensor irreps and their duals }  
We specialize $\Up$ and $\Um$ to the case with ${\bf\Theta}={\bf I}'
\backslash \{m+n-1\}$. 
Consider a one - dimensional irreducible $\Up$ module $V_+=\C v$ 
such that 
\ban 
E_{b\, b+1} v &=& E_{c+1\, c} v =0, \\
K_b v &=& v, \\
K_{m+n} v &=& q^{-k} v, \\
k\in\Z_+, & & b, \ c \in {\bf I}', \ c<m+n-1,  
\nan 
and denote the associated representation by $\phi$. 
Define 
\ban 
{\overline{\cal O}}_k &=& \left\{\zeta\in 
V_+ \otimes G_q^{\pi} \, 
|\, (p\circ\zeta)(x) = \phi(S(p)) \zeta(x), \  
 \forall x\in\Uq, \ p\in\Up\right\}. 
\nan   
Direct calculations can show that 
\ba {\overline{\cal O}}_k &=& \bigoplus_{|L|=k}\C v\otimes Z^L, \na 
where $Z^L$ is defined by (\ref{Z}).   Then ${\overline{\cal O}}_k$ 
gives rise to the rank $k$ skew supersymmetric tensor irrep of $\Uq$, 
with the highest weight 
\ban 
\lambda &=& \left\{\begin{array}{l l} 
            \sum_{i=1}^k \epsilon_i, & k\le m, \\
            \sum_{i=1}^m \epsilon_i + (k-m)\epsilon_{m+1}, & k>m. 
            \end{array}\right. 
\nan

Now let $V_-=\C w$ be a one dimensional irreducible $\Um$ module 
such that 
\ban 
E_{c\, c+1} v &=& E_{b+1\, b} v =0, \\
K_b v &=& v, \\
K_{m+n} v &=& q^k v, \\
k\in\Z_+, & & b, \ c \in {\bf I}', \ c<m+n-1,
\nan  
and denote the corresponding irrep by $\psi$.  Define 
\ban
{\cal O}_k &=& \left\{\zeta\in
V_- \otimes G_q^{\bar\pi}  \,    
|\, (p\circ \zeta)(x) = \psi(S(p)) \zeta(x),\ 
 \forall x\in\Uq, \ p\in\Um\right\}.
\nan
Then 
\ba {\cal O}_k &=& \bigoplus_{|L|=k}\C w\otimes  {\overline Z}^L.\na
This time ${\cal O}_k$ yields an irrep with highest weight  
\ban \lambda&=& -k\epsilon_{m+n},  \nan 
which is  dual to the rank $k$ 
skew supersymmetric tensor irrep.

 
\end{document}